\documentclass{amsart}
\usepackage{graphicx}

\def\lput(#1,#2)#3{\put(#1,#2){\hbox to 0pt{\hss{#3}}}}
\def\allowhyphens{\penalty10000 \hskip0pt}
\def\dash{\allowhyphens\discretionary{-}{}{-}\allowhyphens}

 \title{Vladimir Igorevich Arnold}
\author{ Oleg Karpenkov}
\address{Oleg Karpenkov \\
Institut f.\ Geometrie, TU Graz \\
Kopernikusgasse 24\\
8010 Graz.}

\begin{document}

\maketitle

\hangindent-.5\columnwidth \hangafter=-19
\unitlength.01\columnwidth
\begin{picture}(.01,.01)\lput(96,-56){
        \includegraphics[height=.57\columnwidth]{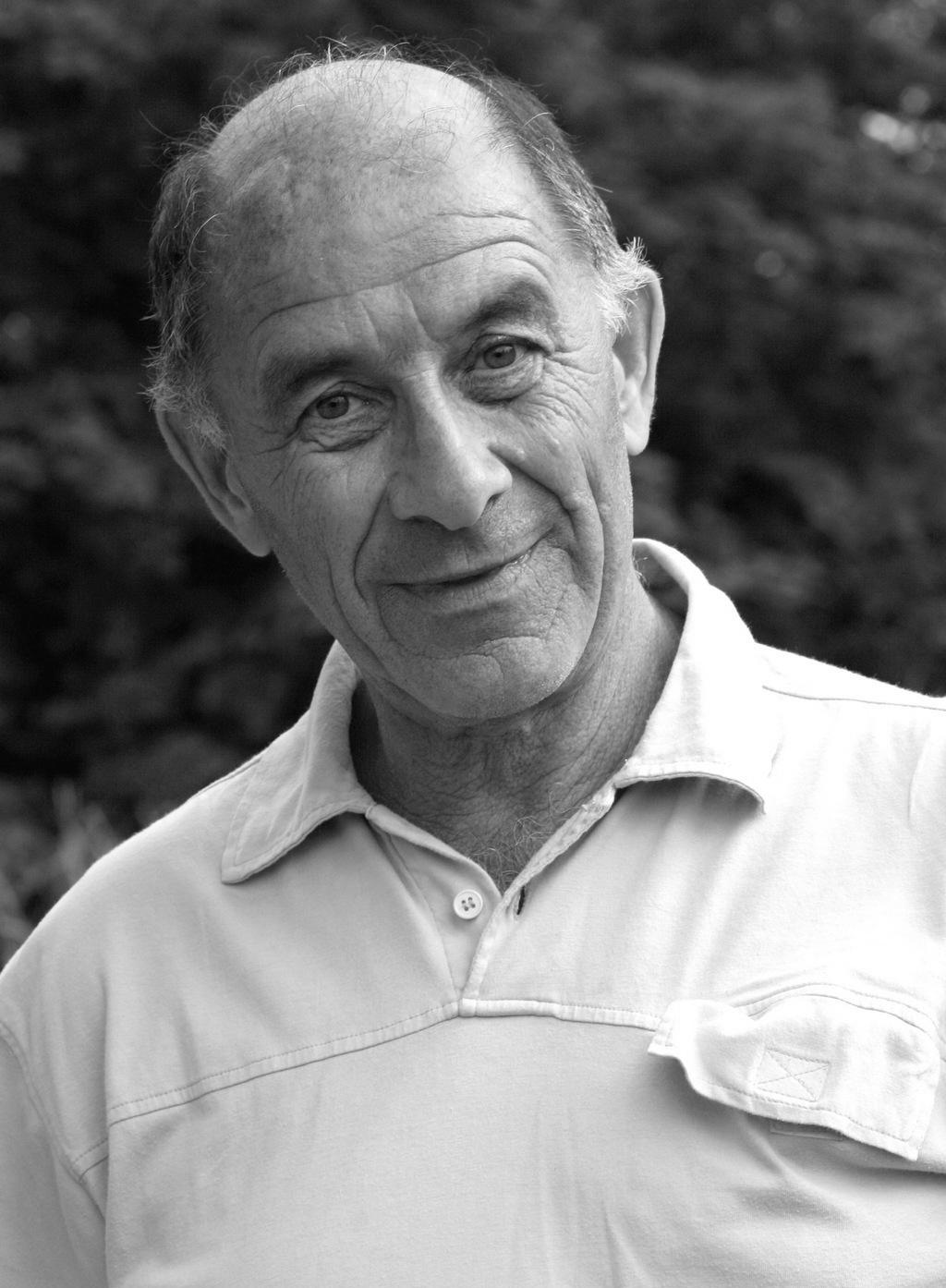}}
\end{picture}%
On June 3, 2010, the
outstanding mathematician Vladimir Igorevich Arnold died in
Paris at the age of 72. He was an
international member of the French Academy of Science since 1983,
an international member of the U.S.~National Academy of Sciences
since 1984, a corresponding member\slash member of the Russian Academy
of Sciences since 1986/1990, a member of the Academy of Arts and
Sciences of the United States since 1987, a member of the London
Royal Society since 1988, an international member of Accademia dei
Lincei in Rome since 1988, and a member of the Russian Academy of
Natural Sciences since 1991. From 1996, V.~I.~Arnold was
President of the Moscow Mathematical Society.

Vladimir Arnold started his mathematical life with finding a
solution of the 13\dash th Hilbert problem at the age of twenty. For this
result he received the Prize of the Moscow Mathematical Society in
1958. He had a great influence on many branches of mathematics
and its applications. He is one of the co\dash authors of the
theory of small denominators (KAM\dash theory), developed a theory of
Lagrangian singularities, and he laid down the basics of symplectic topology.
Vladimir Arnold has influenced differential equations and partial
differential equations, singularity theory, topology, theory of
braids, real algebraic geometry, magneto\dash hydrodynamics, the theory of
multidimensional continued fractions, finite projective geometry,
and combinatorics.

Together with his scientific advisor A.~N.~Kolmogorov, V.~I.~Arnold
received the Lenin prize of the Soviet Government in 1965, for
work on celestial mechanics. In 1982 he received the
Crafoord Prize of the Swedish Academy of Sciences (with
L.~Nirenberg) for his theory of nonlinear differential equations. In
1994 he was awarded the Technion's Harvey prize for a {\it basic
contribution to the stability theory of Dynamical Systems, his
pioneering work on singularity theory and seminal contributions to
analysis and geometry}. Further, in 2001, he received the Dannie
Heineman Prize for Mathematical Physics for {\it his fundamental
contributions to our understanding of dynamics and of
singularities of maps with profound consequences for mechanics,
astrophysics, statistical mechanics, hydrodynamics and optics} and
the Wolf Prize in Mathematics. In 2007 he was awarded the State
Prize of the Russian Federation in science and technology for {\it
outstanding success in mathematics}. Together with
L.~D.~Faddeev, V.~I.~Arnold in 2008 received the Shaw Prize {\it
for their contributions to mathematical physics}.

\vspace{2mm}

{\bf Supervising Ph.D.\ theses.} I thank my lucky stars that I had
an opportunity to be supervised by Vladimir Arnold and to work
with him for almost ten years. Despite all the academic honors and
rewards, Vladimir Arnold always kept in touch with students.
Anyone was able to become a student of Vladimir Arnold, I did not
remember that he refused to work with someone. However, only a few
students continue. It was not complicated to become his student,
but it took an effort to stay his student and later to write a
Ph.D.\ thesis under his supervision.

The first task a student receives from Vladimir Arnold is to
solve the exercises of the so\dash called {\it mathematical trivium}.
This is a list of 100 selected exercises, that any
master student should be able to solve. In addition, Vladimir
Arnold said that any of these problems should take at most 5
minutes. They cover many branches of mathematics, and each
of them is dedicated to some bright idea in the corresponding
area. Let us formulate one of them.

\begin{quote}
\it Problem 2: Compute
$\displaystyle
\lim_{x\to 0}\
\frac{\sin\tan x-\tan \sin x}{\arcsin\arctan
x- \arctan\arcsin x}
$.
\end{quote}

This famous exercise has two essentially different solutions. One
is straightforward: one should apply the rule of de
L'H\^opital. However, the calculation of
derivatives by hand can take several hours or even days depending
on the technical ability of the student. The second solution is
very elegant. I will skip it in order to give the reader the
opportunity to discover the
idea and to enjoy the process of solving the exercise. Usually
each exercise takes much more than just 5 minutes to solve. In
my case it was only 10 exercises per week. Vladimir Arnold checked
the exercises only in written form. If there is even a very small
non\dash clear or non\dash accurate statement in the solution, the exercise
is returned unaccepted with partially humorous remarks. As I heard,
there is no one who successfully solved all (or even 90 of 100) exercises
-- it is almost impossible. While probably
many of these exercises were collected by Vladimir Arnold from
folklore, several definitely belong to himself.

When the majority of exercises is done, Vladimir Arnold has collected
evidence that he works with a strong and diligent student. At that
time he already has some information on the abilities and mathematical
priorities of his new follower. Nevertheless he never proposed
certain problems to his students, he told them that they can
choose any problem he (Arnold) is currently interested in.
When I asked him for a problem for my Ph.D., he
told me that {\it to choose a problem is as personal as
to choose a wife}, so to ask this question to someone else is just useless.
Each year he announced many new problems and new directions that
appeared in mathematics, many of which
are collected in the book ``Arnold's Problems'' ({\it
Springer Verlag, Berlin; PHASIS, Moscow, 2004}), where they are
supplemented by comments on their state at that time by
mathematicians who personally contributed to them.
The half\dash life of such a problem turned out to be
approximately 8 years. Many are still open
or only partially solved.

As a rule, Vladimir Arnold
worked in a certain area for several years. He did there as much as he
could, formulating problems and conjectures for his students
and colleagues for further study. Then he usually switched to
another area, giving to his followers the opportunity to enjoy the beauty of
mathematics. Such a system works perfectly, since
students do not compete with other former students, and the knowledge
of the members of Arnold's school is complementary but not
competitive.

It was usually hard to write a Ph.D.\ thesis under the supervision of
Vladimir Arnold. When you have enough material for the thesis and
finally bring a text to Vladimir Arnold
to read (waiting for his praise and admiration), then within one week he
gives the thesis back to you, together with a review for the PhD
defense, looking at you in a very amicable and gentle\dash ironical way.
So you
are almost ready to go and defend the thesis, but then you suddenly
have a closer look at the returned text and you find
that all the empty spaces on each page and sometimes between the
lines are filled by
Arnold's remarks. The amount of remarks and corrections was
really comparable to the amount of the text itself. The last hope
for easy money is lost when you read his review. First
he describes the area of the thesis and its importance in
science and applications, but then he continues with the
following: {\it unfortunately the student did not understand all
this, he did not solve the original problems, and the text is
completely unreadable.} Finally he concludes with an exposition of
currently unsolved problems in this area for further
investigation and wishes you good luck in the PhD defense.
Of course it is useless even to try to
defend the thesis with such a review. So after this trip back to the earth
you work hard for several weeks and prepare a completely
revised text (that is really an improvement on the
first one). Then you give the second version to Vladimir Arnold,
and he returns within one week with new comments
and a review. In this review he softens some of the original statements,
but usually it is still not sufficient for
a PhD defense. So you should continue to improve your text.
When I had rewritten my thesis 7 times, I finally received a review
with which I could try to defend the thesis. Later I found out that I
was actually lucky: in some cases the number of iterations was 20!
Here I should admit that after each correction cycle the quality
of the text in general improves, and its mathematical worth
increases as well. In addition, the student learns much on how to write.
All of Vladimir Arnold's articles and
books are written in perfect language,
understandable to specialists from different areas: he was aiming
to educate them in the same way as his students.

\vspace{2mm}

{\bf Outlook on Mathematics.} Vladimir Arnold liked to joke that
mathematicians often experience life in the way of a person much
younger than they actually are. They somehow finish their `growth'
at certain age. The younger the age of maturity, the more
brilliant the professional characteristics of a mathematician
usually are. Then he states that he himself stopped growing at the
age of a 12--13 year old boy. This is approximately the age when a
child is wondering about everything he sees. He is trying to play
with everything and to make simple experiments which do not make
sense for an adult, with the purpose to answer thousands of `why'
questions. Usually when children grow older, they are less and
less interested in such things.

Experiments formed an essential part of Vladimir Arnold's
research. He stated that mathematics is an experimental science.
If one is aiming to solve some problem first he should collect
information. Many laws just cannot be seen without an exhaustive
study of examples.

In addition Vladimir Arnold did not like to distinguish the
branches of mathematics, since it is based on bright ideas that
appear in theorems of different areas of mathematics. This can be
easily seen by means of the idea of Euclid's algorithm and its realization
in terms of ordinary continued fractions. The algorithm was
invented in number theory to find the greater common divisor
$\gcd(p,q)$ of
two integers $p$ and $q$, and works as follows.
The first step is to find integers $a_0$ and $r_1$
with $q>r_1\ge 0$ such that
    $$ p=a_0q+r_1. $$
 Suppose we have completed $k-1$ steps, resulting in integers
$a_{k-2}$ and $r_{k-1}$. Then we find $a_{k-1}$ and $r_{k}$ where
$r_{k-1}>r_k\ge 0$ such that
    $$
    r_{k-2}=a_{k-1}r_{k-1}+r_k.
    $$
The algorithm stops after $n+1$ steps, and we have $\gcd(p,q)=r_{n}$.
The key point of the algorithm is computation of the numbers
$a_0,\ldots, a_n$. It is interesting to note that these numbers
form an ordinary continued fraction for $p/q$. Indeed,
    $$
    \frac{p}{q}=a_0+\frac{1}{\displaystyle a_1+\frac{1}{\displaystyle
    a_2+\frac{1}{\displaystyle\ \ddots\ +\frac{\displaystyle 1}{a_n}}}}.
    $$
 To calculate $a_0,\ldots, a_n$ for a positive rational number
${p}/{q}$, we only employ the two basic
operations $s \to s-1$ and $s\to s^{-1}$:
    $$
    \textstyle
\frac{p}{q}=1{+}\Big(\frac{p}{q}{-}1\Big)=\ldots=a_0{+}\Big(\frac{p}{q}{-}\left\lfloor
\frac{p}{q}\Big\rfloor\right)=
a_0{+}\frac{1}{\frac{1}{\frac{p}{q}{-}\left\lfloor
\frac{p}{q}\right\rfloor }}=
a_0{+}\frac{1}{1{+}\left(\frac{1}{\frac{p}{q}{-}\left\lfloor
\frac{p}{q}\right\rfloor}{-}1\right)}=\\
\ldots
$$
Here $\lfloor x \rfloor$ denotes the greatest integer
$\le x$. Similar operations appear in
different problems, and interpretations of Euclid's algorithm are used
in Anosov systems in dynamic theory, in classifications of
rational knots and 3\dash dimensional manifolds in topology,
Hirzebruch\dash Jung theory of singularity resolutions. Simplest
continued fractions appear in the spiral patterns of sunflowers
and cones of trees. This is also due to the fact that the
corresponding biological laws somehow refer to the mentioned
operations.

In the opinion of Vladimir Arnold, the difference between mathematics
and physics is only in the total cost of experiments. In physics,
they may cost millions like in the case of the Large Hadron
Collider, while in mathematics they cost almost nothing: one can
perform them with a pen and paper. In the past centuries
certain areas of mathematics, mechanics, and astronomy were close
and essentially dependent on each other.

Vladimir Arnold payed extra attention to the history of science and
knew many details from the original classical texts of
the great mathematicians of the past centuries: L.~Euler,
H.~Poincair\'e, F.~Klein, etc. He liked to know by whom
a theorem was discovered first. The Arnold principle
states that {\it theorems, theories, principles, etc.\ in science
are usually not named after their authors}. While saying this he
always added that {\it Arnold's principle applies to itself}.
He explained to us that there are two ways to avoid
wrong attributions: One should either tell all
discoveries fairly to everyone, or to keep them secret up to the
end. In both cases there is no opportunity for another party
to steal the results. He mentioned that he prefers and follows
only the first strategy and taught us not to be lazy to give talks
and to make them in such a way that everyone can understand them,
even if he works in another field of mathematics.

\vspace{2mm}

{\bf Arnold's seminars.} Vladimir Igorevich was running two
permanent seminars: one in Moscow and one in Paris. I attended the
Paris seminar only a small number of times, but I was regularly
participating in the Moscow seminar which took place in room
14\dash 14 of the Main Building of Moscow State University. This
building is situated on top of a spectacular green hill (Vorobyevi
Gori) to the South\dash West of the center of Moscow near the
Moscow river. It is one of seven Stalinist neoclassic\dash style
skyscrapers built in 1953. It has 36 floors, the Mathematical
Department being located at floors 12--16, with a magnificent view
of the city. A seminar usually started on Tuesdays at around 4
p.m.\ and officially lasted for one and a half hours. In former
times it could continue for two or three more hours, depending on
the subject. When I started to attend the seminar in 1998, it
lasted usually no more than two hours, since many of its
participants, including Arnold himself, were involved in the
meetings of the Moscow Mathematical Society scheduled two hours
after the seminar.

The first talk of each winter semester was always given by Vladimir
Arnold himself. This was the most important event in the life of
the seminar. In this talk he discussed new
problems for the seminar that he had collected recently. The list
usually contained 10--20 problems on 4 or 5 different subjects.
Some of them were rather elementary and perfectly suitable for
students starting their research, while other problems were dealing
with more complicated subjects that sometimes even gave a new
sight on classical problems and conjectures. These problems form
the core of the volume ``Arnold's Problems'' mentioned above.

Apart from the problems, Vladimir Arnold usually distributed
`tasks' to his former students. He brought recent
papers that attracted his
attention and asked the members of the seminar if someone is
interested to read some of them and give their opinion. In the
majority of cases the papers were interesting, and we spent one seminar or
half of a seminar to listen to the participants who had read
the papers.

The remaining seminar talks were given by his students and
colleagues who presented their new results, and in part by distinguished
colleagues who gave survey talks. This seminar did not have a certain
fixed subject. Vladimir Arnold himself contributed to many
branches of mathematics, including applied mathematics and mechanics. So one
could expect talks on current problems all over mathematics and
its applications. It was felt to be an honor to speak
at the seminar and to receive
comments of Vladimir Arnold and his former students, many of whom
are now holding professor positions at Moscow State University or
some other institutes or universities in Moscow. Almost always the
feedback received from the audience led to new insights, additional
references, and ideas for future work.

It was really a
peculiar experience for a speaker to have Vladimir Igorevich in the audience.
He would interrupt any speaker, especially in his seminar.
When the talk started without taking care of all important definitions
and notions, he would ask for all the necessary details and himself
explain the importance of these details to the speaker, mentioning
examples when different understandings or definitions of certain
notions lead to completely different results. Such a speaker would be lucky
to give the main definitions and the main result in one and a half hours.
In some sense, of course, this is the minimal amount of
information necessary for the audience. So at the end of the talk
it was completely clear what was the result and what are the open questions
(which are often not mentioned by the
speakers, but are extremely interesting for the audience). While
listening to such a presentation, Vladimir Arnold did not distinguish
between a young unknown student and a famous academician. Of
course not all people would appreciate such an attitude
and from time to time they got quite angry (Vladimir Igorevich
did not pay attention to their complaints). All this made the talks extremely
useful for the audience and in many cases for the speaker himself.
Sometimes, when the speaker ignored his remarks, Arnold
would start to correct his students' papers or think about his own research.
Thus the speaker was in some sense
`punished': his talk had become an ordinary talk, with the
majority of the audience not paying attention, either
switching to their own businesses, or diligently
continuing looking at the blackboard with eyes wide open.

\vspace{2mm}

{\bf On abstraction and teaching.} Vladimir Arnold was against
what he called `algebraization' of mathematics which is very
popular in our days among many people whom he would refer to as
the {\it devils of algebra}, while he attributed himself as the
{\it angel of geometry}. One cannot conclude from this that he was
`against' certain areas of mathematics, but he was against
eliminating the thinking process by formal calculations, which
apparently he thought certain followers of N.~Bourbaki are guilty
of.

Once I was reading a preprint of a school
textbook, and there was a standard question on percentages:
{\it with a monthly interest rate of
1$\%$, by how many percent has the capital increased after one year?}
The solution given by the author was
around 20000$\%$. This is clearly not true, and
you would never write such nonsense if you were thinking while writing the
book. Unfortunately similar situations sometimes present themselves in
mathematical articles where catching mistakes is not as
easy.

Vladimir Igorevich did always care about the future of
mathematics. On the one hand he was a person with great knowledge
which he was trying to transfer to his colleagues, students, or
even persons met by chance. On the other
hand, he clearly understood that the roots of mathematics lie in
education. If we neglect the quality of teaching in schools,
we will get a generation unable to
expand, usefully employ, or even to keep the achievements of our and
previous generations. A recent issue here is the use of computers which
in many cases makes it possible to
press the correct buttons instead of thinking. As a result we get pupils
unable to do simple operations without a computer.
Here are my favorite examples of `wrong' answers given by pupils
which I know from Vladimir Arnold.
The first one is coming from some American test where it turned out that
the majority of schoolchildren did not know to handle fractions: They did
computations like
    $$
    \frac{2}{3}+\frac{1}{2}=\frac{2+1}{3+2}=\frac{3}{5}.
    $$
 The second example is related to the traditional
emphasis on axiomatics in France. Vladimir Arnold asked a pupil with
an excellent record in mathematics: {\it what is $2+3$}? After a
minute of thinking the pupil replied that
    $$
    2+3=3+2,
    $$
since addition is a commutative operation.
Vladimir Arnold was constantly trying to improve
mathematical education, writing many
articles, critical reviews, and giving public talks with the purpose of
educating a new generation of mathematicians.

\vspace{2mm}

{\bf Private life and work.} It is hard to believe how many things
Vladimir Arnold apparently managed to do. The reviewing of theses
and articles of his students and colleagues was just a small part
of his work. In the last years he wrote what must have been up to
100 handwritten pages per day. This includes new articles, books
on mathematics, articles on popular mathematics and recent
problems of education. It is interesting to know that Vladimir
Arnold did not like to work on a computer in his office, but he
preferred to do research with pencil and paper, occupying one of
his favorite outdoor places: sitting under the oldest tree in
Paris or somewhere in the Bois de Vincennes; somewhere in the
rocky region of Trieste, where he took part in archeological
research in the local caves; at his dacha near Moscow in an area
with beautiful lakes, rivers and forests, berries and mushrooms.

I have experienced Vladimir Arnold as a perfect person to discuss any subject.
Communication with him was a powerful stimulus for
research, and yielded inspiration for work even months after
it actually took place. He preferred to meet people in an informal environment,
for instance in his dacha, where the meeting
could last all evening. His wife Eleonora
Aleksandrovna is a hospitable, very friendly and kind person,
who used to prepare delicious dinners for their guests.
After such a dinner
Vladimir Arnold invited his guests for a round trip in the forest, on
foot, on skis, or by bike.
In Paris he usually invited his guests and colleagues for 4--5 hour
excursions and gave bright expositions of history that were full
of stories about the adventures of G.~Casanova, or the guillotine, or
or the poisons of
Catherine de Medici.

In Russia, Vladimir Arnold
spent most of the time at his famous dacha situated in the
academician suburban settlement not far from Moscow. To go there, one takes
the bus from the local train station {\it Perkhushkovo} and alights
at bus stop ``Dacha No.~1'' (Dacha No.~1 is one of
the dachas belonging to Stalin).
In such settlements academicians can meet and talk to each other
privately and informally, and you
could meet leading people of all branches of science.
Doubtless this invitation to multidisciplinary collaboration was
appreciated by many Russian scientists living there.

When Vladimir Arnold stayed in Moscow, his students, friends, and
colleagues had an opportunity to meet another challenge:
It sounds quite simple to run 40
kilometers on cross\dash country skies together with Vladimir Arnold
and his older students in in the forests on
both sides of Moscow river. But the main obstacle to completing this
particular course is that it
starts on one side of the river and
finishes on the other. The river in this place is quite close to
its origin and not very wide, but fast enough not to
freeze in most winters. Approximately halfway,
Vladimir Arnold crossed the river
by swimming, with the majority of his group following him (I
should add at this point that the only thing Vladimir Arnold
was wearing were his shorts). In his best days, Vladimir Arnold was in very
good shape physically: he used to go by bike but stated that mere
bike roads were too simple, he would swim 10km from Luminy to
Marseille, or sail in a canoe for several
days. Of course like every man of his age, he had some health problems.
Still all who knew him were sure that he would
celebrate his 80\dash th, 90\dash th, and 100\dash th anniversaries. He died
unexpectedly. Several days before his death he was full of plans,
new ideas and enthusiasm. Now he has left us. But we still have
his rich intellectual legacy contained in his many
lecture notes, books, articles and students, who were
inspired by the enthusiastic, optimistic, and always young genius which
Vladimir Arnold definitely was.

\bigskip

{\bf Acknowledgement.} The author is grateful for S.~Tretyakova
for the photo of Vladimir Arnold used in this article.

\end{document}